\documentclass[a4paper,11pt]{amsart}

\usepackage{amsfonts, latexsym, amsmath, amssymb, amsthm, amscd}
\usepackage[all]{xy}
\usepackage[english]{babel}
\usepackage[utf8]{inputenc}
\usepackage[T1]{fontenc}
\usepackage{graphicx}
\usepackage{booktabs}
\usepackage{braket}
\usepackage{mathtools}
\usepackage{array, tabularx}
\usepackage{setspace}
\usepackage{textcomp}
\usepackage{tikz}
\usepackage{multirow}
\usepackage{paralist}
\usepackage{comment}
\usepackage{url}
\usepackage{tensor}

\usepackage[hypertexnames=false,
backref=page,
    pdftex,
    pdfpagemode=UseNone,
    breaklinks=true,
    extension=pdf,
    colorlinks=true,
    linkcolor=blue,
    citecolor=blue,
    urlcolor=blue,
]{hyperref}

\usepackage[top=0.5in, bottom=.9in, left=0.9in, right=0.9in]{geometry}

\newcommand{\referenza}{}

\newtheorem{remark}{Remark}
\newenvironment{re}[1]{\renewcommand\theremark{#1}\remark}{\endremark}

\newtheorem*{thm*}{Theorem (\referenza)}

\newtheorem*{cor*}{Corollary \referenza}

\newtheorem*{lem*}{Lemma \referenza}

\newtheorem*{prop*}{Proposition \referenza}

\newtheorem*{conj*}{Conjecture \referenza}

\DeclareMathOperator{\rank}{rank}

\newcommand{\ud}{\mathop{}\!\mathrm{d}}
\newcommand{\p}{\mathop{}\!\mathrm{\partial}}
\newcommand{\pb}{\mathop{}\!\mathrm{\bar{\partial}}}

\newcommand{\ph}{\varphi}

\newcommand{\ou}{\bar{1}}
\newcommand{\od}{\bar{2}}

\newcommand{\phb}{\bar{\varphi}}

\newcommand{\w}{\wedge}
\newcommand{\g}{\mathfrak{g}}
\newcommand{\h}{\mathfrak{h}}

\newcommand{\R}{\mathbb{R}}
\newcommand{\I}{\mathbb{I}}
\newcommand{\Z}{\mathbb{Z}}
\newcommand{\C}{\mathbb{C}}
\newcommand{\del}{\partial}
\newcommand{\delbar}{\bar{\partial}}
\newcommand{\N}{\mathbb{N}}
\newcommand{\tg}{\mathbf{t}}
\numberwithin{equation}{section}

\allowdisplaybreaks

\begin{document}
\title[On The Fr\"olicher spectral sequence of the Iwasawa manifold]{On the Fr\"olicher spectral sequence of the Iwasawa manifold and its small deformations}
\author{Cosimo Flavi}
\begin{abstract}
We determine the successive pages of the Fr\"olicher spectral sequence of the Iwasawa manifold and some of its small deformations, providing new examples and counterexamples on its properties, including  the behaviour under small deformations.
\end{abstract}
\maketitle

\section*{Introduction}
Thanks to a result by J. Stelzig (\cite[Theorem 3]{Ste18}) we obtain the explicit algebraic structures of the complexes of differential forms of the Iwasawa manifold and some of its small deformations. We introduce these examples whit the aim of studying some properties of the Fr\"olicher spectral sequence; in particular, for each deformation analysed, we compute the explicit dimensions of the components of the successive pages. We observe that:
\begin{enumerate}[(i)]
\item the Fr\"olicher spectral sequence of the Iwasawa manifold degenerates at the second page, but there are some deformations for which it degenerates at the first page;
\item for some deformations the dimensions of the components of the second page of the Fr\"olicher spectral sequence can be in general higher or lower than those of the corresponding components for the Iwasawa manifold;
\item for each page of the Fr\"olicher spectral sequence and for every $p,q\in\Z$ the component of type $(p,q)$ has the same dimension as the component of type $(3-p,3-q)$.
\end{enumerate}

Fact (i) was already known, by considering the Hodge numbers of the Iwasawa manifold and its small deformations determined by I. Nakamura in \cite[p. 96]{Nak75}, in which it is shown that the degeneration at the second page is not a stable property under small deformations. The fact (ii) shows instead that the dimensions of the various components do not vary in general either upper or lower semi-continuosly under small deformations. 

These two considerations were already provided by M. Maschio in \cite{Mas18}, where he considered the Nakamura manifold, which is an example of solvmanifold, and a family of its small deformations. In particular, in that case the behaviour of the Fr\"olicher spectral sequence is different; in fact, it degenerates at the second page for the Nakumara manifold, while it degenerates at higher pages for every other deformations.

Finally, fact (iii) is a direct consequence of a generalization of the Serre duality, recently proved by A. Milojevic in \cite{Mil19}.

\section{Preliminaries}
A double complex over a field $K$ is a triad $(A,\del_1,\del_2)$, where $A=\bigoplus_{p,q\in\Z}A^{p,q}$ is a bigraded $K$-vector space and $\del_1$, $\del_2$ are two endomorphisms of bidegree $(1,0)$ and $(0,1)$, i.e. such that for each $p,q\in\Z$
\[
\del_1(A^{p,q})\subseteq A^{p+1,q},\quad \del_2(A^{p,q})\subseteq A^{p,q+1},
\]
satisfying the \emph{boundary conditions} $\del_1^2=0$, $\del_2^2=0$ and also
 the anti-commutativity property $\del_1\del_2+\del_2\del_1=0$. From now by the term \emph{complex} we will refer to a double complex over a field $K$. 

A result by J. Stelzig (\cite[Theorem 3]{Ste18}) guarantees that every bounded complex can be decomposed in an essentially unique way as direct sum of \emph{indecomposable} complexes, i.e. those complexes that cannot be seen as direct sum of two other nontrivial complexes. This class is composed by \emph{squares}, defined as those complexes whose representation is of the type
\[
\xymatrix{
A^{p,q+1} \ar[r]^{\del_1} & A^{p+1,q+1} \\
A^{p,q} \ar[u]^{\del_2} \ar[r]^{\del_1} & A^{p+1,q} \ar[u]^{\del_2}
},
\]
and \emph{zigzags}, given by
\begin{gather*}
\xymatrix{
A^{p,q}
},\quad
\xymatrix{
A^{p,q} \ar[r]^{\del_1} & A^{p+1,q}
},\quad
\xymatrix{
A^{p,q+1} \\
A^{p,q} \ar[u]^{\del_2}
},\quad
\xymatrix{
A^{p,q+1} & \\
A^{p,q} \ar[u]^{\del_2} \ar[r]^{\del_1} & A^{p+1,q}
}, \\
\xymatrix{
A^{p,q+1} \ar[r]^{\del_1} & A^{p+1,q+1} \\
& A^{p+1,q} \ar[u]^{\del_2}
},\quad
\xymatrix{
A^{p,q+1} \ar[r]^{\del_1} & A^{p+1,q+1}  \\
& A^{p+1,q} \ar[u]^{\del_2} \ar[r]^{\del_1} & A^{p+2,q+1}
},\quad \dots
\end{gather*}
where in both cases all the drawn components represent one-dimensional vector spaces and all drawn arrows represent isomorphisms. The number of non-zero components is said \emph{lenght} of the zigzag.

Let us consider now the complex of complex-valued differential forms on a complex manifold $X$, denoted by $\big(\Omega_\C(X),\del,\delbar\big)$. If $X=\Gamma\backslash G$ is a \emph{nilmanifold}, i.e. a compact quotient of a connected and simply connected nilpotent Lie group $G$ by a co-compact discrete subgroup $\Gamma$, endowed with a left-invariant complex structure then, under some appropriate hypotheses, all cohomological informations can be obtained from the subcomplex of left-invariant forms, namely differential forms whose pullback to $G$ is invariant by left-translations, denoted by $\big(\Omega_\C^\text{inv}(X),\del,\delbar\big)$ (see \cite[Chapter 3]{Ang14} and references therein), which is isomorphic to the complex $\big(\bigwedge^{\bullet,\bullet} \g_\C^*,\del,\delbar\big)$, where $\g_\C\coloneqq\g\otimes\C$ represents the complexification of the Lie algebra $\g$ naturally associated to $G$. Using the results obtained by L. A. Cordero, M. Fernández, L. Ugarte and A. Gray in \cite[Theorem 1, Theorem 3]{CFUG97} this last complex also permits us to compute easily the dimensions of the successive pages of the Fr\"olicher spectral sequence ${\{(E_r^{\bullet,\bullet},\ud_r)\}}_{r\in\N}$.

A classical example of nilmanifold endowed with a left-invariant complex structure is the \emph{Iwasawa manifold}, defined as the quotient 
\[
\mathbb{I}\coloneqq\mathbb{H}\big(3;\Z [i]\big)\big\backslash\mathbb{H}(3;\C),
\]
where $\mathbb{H}(3;\C)$ is the $3$-dimensional \emph{Heisenberg group} over $\C$, defined by
\[
\mathbb{H}(3;\C)\coloneqq\Set{
\begin{pmatrix}
1 & z_1 & z_3 \\
0 & 1 & z_2 \\
0 & 0 & 1
\end{pmatrix} \in \mathrm{GL}(n,\C)
| z_1,z_2,z_3\in\C}
\]
and $\mathbb{H}(3,\Z[i])$ corresponds to the lattice
\[
\mathbb{H}(3,\Z[i])\coloneqq\mathbb{H}(3,\C)\cap \mathrm{GL}(3,\Z[i]).
\] 
I. Nakamura determined in \cite[pp. 94-96]{Nak75} a family of small deformations $\{X_{\mathbf{t}}=(\I,J_{\mathbf{t}})\}_{\mathbf{t}\in\Delta(0,\varepsilon)}$ depending on six parameters 
\[
\mathbf{t}=(t_{11},t_{12},t_{21},t_{22},t_{31},t_{32})\in\Delta(0,\varepsilon),
\]
with
\[
\Delta(0,\varepsilon)\coloneqq\Set{\mathbf{s}\in\C^6|\lvert \mathbf{s}\rvert <\varepsilon},
\] 
where $\varepsilon>0$ is small enough.

Furthermore, he provided a co-frame $(\ph_\tg^1,\ph_\tg^2,\ph_\tg^3)$ of left-invariant $(1,0)$-forms on $X_\tg$ satisfying the structure equations given by
\[
\begin{cases}
\ud\ph_\tg^1=0\\
\ud\ph_\tg^2=0\\
\ud\ph_\tg^3=\sigma_{12}\,\ph_\tg^1\w\ph_\tg^2+\sigma_{1\ou}\,\ph_\tg^1\w\phb_\tg^1+\sigma_{1\od}\,\ph_\tg^1\w\phb_\tg^2+\sigma_{2\ou}\,\ph_\tg^2\w\phb_\tg^1+\sigma_{2\od}\,\ph_\tg^2\w\phb_\tg^2,
\end{cases},
\]
where $\sigma_{12}, \sigma_{1\ou}, \sigma_{2\od}, \sigma_{2\ou}, \sigma_{1\od}\in\C$ are parameters depending on $\tg$. 
Defining the matrix $S$ as
\[
S=
\begin{pmatrix}
\overline{\sigma_{1\ou}} & \overline{\sigma_{2\od}} &\overline{\sigma_{1\od}} &\overline{\sigma_{2\ou}}  \\
\sigma_{1\ou} & \sigma_{2\od} & \sigma_{2\ou} & \sigma_{1\od}
\end{pmatrix},
\]
the small deformations of the Iwasawa manifold can be classified into three classes, based on their Hodge numbers, further subdivided in other subclasses considering their Bott-Chern-numbers:
\begin{enumerate}[(i)]
\item $t_{11}=t_{12}=t_{21}=t_{22}=0$;
\item $D(\tg)=0$ and $(t_{11},t_{12},t_{21},t_{22})\neq(0,0,0,0)$:
\begin{enumerate}[(i)]
\item[(ii.a)] $D(\tg)=0$ and $\rank S=1$;
\item[(ii.b)] $D(\tg)=0$ and $\rank S=2$;
\end{enumerate}
\item[(iii)] $D(\tg)\neq 0$:
\begin{enumerate}[(i)]
\item[(iii.a)] $D(\tg)\neq 0$ and $\rank S=1$;
\item[(iii.b)] $D(\tg)\neq 0$ and $\rank S=2$.
\end{enumerate}
\end{enumerate}

Thanks to some results contained in \cite{Rol09}, \cite{Ang13} and the references therein, we can make some considerations about the de Rham, Dolbeault, Aeppli and Bott-Chern cohomologies, in particular for each $\tg\in\Delta(0,\varepsilon)$ we get the isomorphisms
\begin{gather*}
H_{dR}^\bullet(\h_,\R)\simeq H^\bullet_{dR}(X_\tg,\R),\quad H_{\delbar}^{\bullet,\bullet}(\h_\C
)\simeq H^{\bullet,\bullet}_{\delbar}(X_\tg),\\
H_{BC}^{\bullet,\bullet}(\h_\C)\simeq H^{\bullet,\bullet}_{BC}(X_\tg),\quad H_{A}^{\bullet,\bullet}(\h_\C)\simeq H^{\bullet,\bullet}_{A}(X_\tg),
\end{gather*}
where $\h$ is the Lie Algebra naturally associated to $\mathbb{H}(3,\C)$ and $H_*^{\bullet,\bullet}(\h_\C)$ represents the corresponding cohomology of the subcomplex of left-invariant forms; therefore, to compute dimensions of the cohomologies, it is sufficient to analyze the complex of left-invariant forms. For this reason, in the next section we will restrict our attention on this last complex.

\section{Double complex of left-invariant forms on the Iwasawa manifold and its small deformations}
We are now going to analyze some explicit examples of complex of differential forms; in particular we will give an explicit graphic representation, as in the Stelzig process, of the structure of the complex related to the Iwasawa manifold and some of its small deformations. Once we describe the required decomposition as a direct sum of squares and zigzags, it will be easy to obtain the dimensions of the successive pages of the Fr\"olicher spectral sequence.

Given any arbitrary bounded complex, the decomposition described by J. Stelzig is based on the costruction of an ascending filtration, which further splits into simpler pieces, thus allowing to identify squares and zigzags appearing in the complex. The procedure works as follows: we determine for each $k\in\Z$ a subcomplex generated by all components of total degree at most $k$ and then a complement of the subcomplex generated by the components of total degree at most $k-1$ contained in it. This permits to obtain for every $k\in\Z$ a simpler subcomplex composed by components of degree $k,k+1$ and $k+2$, whose direct sum corresponds to complex considered. As a result it is easier to determine the various squares and zigzags contained in every subcomplex to get the desired decomposition. 

In the case of the complex of left-invariant differential forms of the Iwasawa manifold and its small deformations the procedure consists in five steps, based respectively on the computation of the images of the left-invariant differential forms of total degree between $2$ and $5$ by the operators $\p$ and $\pb$.

\begin{figure}[ht]
\begin{center}
a)\quad
\begin{tikzpicture}[scale=1.5]

\draw[help lines,black] (0,0) grid (4,4);

\draw [line width=1.5,->,red] (3/2,1/2) -- (5/2-0.03,1/2);
\draw [line width=1.5,->,red] (1/2,3/2) -- (1/2,5/2-0.03);

\draw [line width=1.5,->,green] (7/4,5/4) -- (9/4-0.03,5/4);
\draw [line width=1.5,->,green] (7/4,3/2) -- (9/4-0.03,3/2);
\draw [line width=1.5,->,green] (5/4,7/4) -- (5/4,9/4-0.03);
\draw [line width=1.5,->,green] (3/2,7/4) -- (3/2,9/4-0.03);
\draw [line width=1.5,->,green] (7/4,7/4) -- (9/4-0.03,7/4);
\draw [line width=1.5,->,green] (7/4,7/4) -- (7/4,9/4-0.03);
\draw [line width=1.5,->,green] (9/4,7/4) -- (9/4,9/4-0.03);
\draw [line width=1.5,->,green] (7/4,9/4) -- (9/4-0.03,9/4);

\draw [line width=1.5,->,blue] (7/4,5/2) -- (9/4-0.03,5/2);
\draw [line width=1.5,->,blue] (7/4,11/4) -- (9/4-0.03,11/4);
\draw [line width=1.5,->,blue] (5/2,7/4) -- (5/2,9/4-0.03);
\draw [line width=1.5,->,blue] (11/4,7/4) -- (11/4,9/4-0.03);

\draw [line width=1.5,->,orange] (7/2,3/2) -- (7/2,5/2-0.03);
\draw [line width=1.5,->,orange] (3/2,7/2) -- (5/2-0.03,7/2);

\draw [fill] (1/2,1/2) circle [radius=0.03];

\draw [fill=red] (5/4,1/4) circle [radius=0.03];
\draw [fill=red] (3/2,1/2) circle [radius=0.03];
\draw [fill=red] (7/4,3/4) circle [radius=0.03];

\draw [fill=green] (9/4,1/4) circle [radius=0.03];
\draw [fill=red] (5/2,1/2) circle [radius=0.03];
\draw [fill=green] (11/4,3/4) circle [radius=0.03];

\draw [fill=blue] (7/2,1/2) circle [radius=0.03];

\draw [fill=red] (1/4,5/4) circle [radius=0.03];
\draw [fill=red] (1/2,3/2) circle [radius=0.03];
\draw [fill=red] (3/4,7/4) circle [radius=0.03];

\draw [fill=green] (5/4,5/4) circle [radius=0.03];
\draw [fill=green] (3/2,5/4) circle [radius=0.03];
\draw [fill=green] (7/4,5/4) circle [radius=0.03];
\draw [fill=green] (5/4,3/2) circle [radius=0.03];
\draw [fill=green] (3/2,3/2) circle [radius=0.03];
\draw [fill=green] (7/4,3/2) circle [radius=0.03];
\draw [fill=green] (5/4,7/4) circle [radius=0.03];
\draw [fill=green] (3/2,7/4) circle [radius=0.03];
\draw [fill=green] (7/4,7/4) circle [radius=0.03];

\draw [fill=green] (9/4,5/4) circle [radius=0.03];
\draw [fill=blue] (5/2,5/4) circle [radius=0.03];
\draw [fill=blue] (11/4,5/4) circle [radius=0.03];
\draw [fill=green] (9/4,3/2) circle [radius=0.03];
\draw [fill=blue] (5/2,3/2) circle [radius=0.03];
\draw [fill=blue] (11/4,3/2) circle [radius=0.03];
\draw [fill=green] (9/4,7/4) circle [radius=0.03];
\draw [fill=blue] (5/2,7/4) circle [radius=0.03];
\draw [fill=blue] (11/4,7/4) circle [radius=0.03];

\draw [fill=orange] (13/4,5/4) circle [radius=0.03];
\draw [fill=orange] (7/2,3/2) circle [radius=0.03];
\draw [fill=orange] (15/4,7/4) circle [radius=0.03];

\draw [fill=green] (1/4,9/4) circle [radius=0.03];
\draw [fill=red] (1/2,5/2) circle [radius=0.03];
\draw [fill=green] (3/4,11/4) circle [radius=0.03];

\draw [fill=green] (5/4,9/4) circle [radius=0.03];
\draw [fill=green] (3/2,9/4) circle [radius=0.03];
\draw [fill=green] (7/4,9/4) circle [radius=0.03];
\draw [fill=blue] (5/4,5/2) circle [radius=0.03];
\draw [fill=blue] (3/2,5/2) circle [radius=0.03];
\draw [fill=blue] (7/4,5/2) circle [radius=0.03];
\draw [fill=blue] (5/4,11/4) circle [radius=0.03];
\draw [fill=blue] (3/2,11/4) circle [radius=0.03];
\draw [fill=blue] (7/4,11/4) circle [radius=0.03];

\draw [fill=green] (9/4,9/4) circle [radius=0.03];
\draw [fill=blue] (5/2,9/4) circle [radius=0.03];
\draw [fill=blue] (11/4,9/4) circle [radius=0.03];
\draw [fill=blue] (9/4,5/2) circle [radius=0.03];
\draw [fill=orange] (5/2,5/2) circle [radius=0.03];
\draw [fill=orange] (11/4,5/2) circle [radius=0.03];
\draw [fill=blue] (9/4,11/4) circle [radius=0.03];
\draw [fill=orange] (5/2,11/4) circle [radius=0.03];
\draw [fill=orange] (11/4,11/4) circle [radius=0.03];

\draw [fill=magenta] (13/4,9/4) circle [radius=0.03];
\draw [fill=orange] (7/2,5/2) circle [radius=0.03];
\draw [fill=magenta] (15/4,11/4) circle [radius=0.03];

\draw [fill=blue] (1/2,7/2) circle [radius=0.03];

\draw [fill=orange] (5/4,13/4) circle [radius=0.03];
\draw [fill=orange] (3/2,7/2) circle [radius=0.03];
\draw [fill=orange] (7/4,15/4) circle [radius=0.03];

\draw [fill=magenta] (9/4,13/4) circle [radius=0.03];
\draw [fill=orange] (5/2,7/2) circle [radius=0.03];
\draw [fill=magenta] (11/4,15/4) circle [radius=0.03];

\draw [fill] (7/2,7/2) circle [radius=0.03];

\end{tikzpicture}\qquad b)\quad
\begin{tikzpicture}[scale=1.5]

\draw[help lines,black] (0,0) grid (4,4);

\draw [line width=1.5,->,red] (5/4,1/4) -- (9/4-0.03,1/4);
\draw [line width=1.5,->,red] (1/4,5/4) -- (1/4,9/4-0.03);
\draw [line width=1.5,->,red] (1/4,5/4) -- (5/4-0.03,5/4);
\draw [line width=1.5,->,red] (5/4,1/4) -- (5/4,5/4-0.03);

\draw [line width=1.5,->,green] (1/2,5/2) -- (5/4-0.03,5/2);
\draw [line width=1.5,->,green] (7/4,5/4) -- (5/2-0.03,5/4);
\draw [line width=1.5,->,green] (7/4,3/2) -- (9/4-0.03,3/2);
\draw [line width=1.5,->,green] (5/4,7/4) -- (5/4,5/2-0.03);
\draw [line width=1.5,->,green] (5/2,1/2) -- (5/2,5/4-0.03);
\draw [line width=1.5,->,green] (3/2,7/4) -- (3/2,9/4-0.03);
\draw [line width=1.5,->,green] (7/4,7/4) -- (9/4-0.03,7/4);
\draw [line width=1.5,->,green] (7/4,7/4) -- (7/4,9/4-0.03);
\draw [line width=1.5,->,green] (9/4,7/4) -- (9/4,9/4-0.03);
\draw [line width=1.5,->,green] (7/4,9/4) -- (9/4-0.03,9/4);

\draw [line width=1.5,->,blue] (7/4,5/2) -- (9/4-0.03,5/2);
\draw [line width=1.5,->,blue] (11/4,3/2) -- (7/2-0.03,3/2);
\draw [line width=1.5,->,blue] (3/2,11/4) -- (9/4-0.03,11/4);
\draw [line width=1.5,->,blue] (5/2,7/4) -- (5/2,9/4-0.03);
\draw [line width=1.5,->,blue] (3/2,11/4) -- (3/2,7/2-0.03);
\draw [line width=1.5,->,blue] (11/4,3/2) -- (11/4,9/4-0.03);

\draw [line width=1.5,->,orange] (15/4,7/4) -- (15/4,11/4-0.03);
\draw [line width=1.5,->,orange] (7/4,15/4) -- (11/4-0.03,15/4);
\draw [line width=1.5,->,orange] (11/4,11/4) -- (15/4-0.03,11/4);
\draw [line width=1.5,->,orange] (11/4,11/4) -- (11/4,15/4-0.03);


\draw [fill] (1/2,1/2) circle [radius=0.03];

\draw [fill=red] (5/4,1/4) circle [radius=0.03];
\draw [fill=red] (3/2,1/2) circle [radius=0.03];
\draw [fill=red] (7/4,3/4) circle [radius=0.03];

\draw [fill=red] (9/4,1/4) circle [radius=0.03];
\draw [fill=green] (5/2,1/2) circle [radius=0.03];
\draw [fill=green] (11/4,3/4) circle [radius=0.03];

\draw [fill=blue] (7/2,1/2) circle [radius=0.03];


\draw [fill=red] (1/4,5/4) circle [radius=0.03];
\draw [fill=red] (1/2,3/2) circle [radius=0.03];
\draw [fill=red] (3/4,7/4) circle [radius=0.03];

\draw [fill=red] (5/4,5/4) circle [radius=0.03];
\draw [fill=green] (3/2,5/4) circle [radius=0.03];
\draw [fill=green] (7/4,5/4) circle [radius=0.03];
\draw [fill=green] (5/4,3/2) circle [radius=0.03];
\draw [fill=green] (3/2,3/2) circle [radius=0.03];
\draw [fill=green] (7/4,3/2) circle [radius=0.03];
\draw [fill=green] (5/4,7/4) circle [radius=0.03];
\draw [fill=green] (3/2,7/4) circle [radius=0.03];
\draw [fill=green] (7/4,7/4) circle [radius=0.03];

\draw [fill=green] (5/2,5/4) circle [radius=0.03];
\draw [fill=blue] (11/4,5/4) circle [radius=0.03];
\draw [fill=green] (9/4,3/2) circle [radius=0.03];
\draw [fill=blue] (5/2,3/2) circle [radius=0.03];
\draw [fill=blue] (11/4,3/2) circle [radius=0.03];
\draw [fill=green] (9/4,7/4) circle [radius=0.03];
\draw [fill=blue] (5/2,7/4) circle [radius=0.03];

\draw [fill=orange] (13/4,5/4) circle [radius=0.03];
\draw [fill=blue] (7/2,3/2) circle [radius=0.03];
\draw [fill=orange] (15/4,7/4) circle [radius=0.03];


\draw [fill=red] (1/4,9/4) circle [radius=0.03];
\draw [fill=green] (1/2,5/2) circle [radius=0.03];
\draw [fill=green] (3/4,11/4) circle [radius=0.03];

\draw [fill=green] (3/2,9/4) circle [radius=0.03];
\draw [fill=green] (7/4,9/4) circle [radius=0.03];
\draw [fill=green] (5/4,5/2) circle [radius=0.03];
\draw [fill=blue] (3/2,5/2) circle [radius=0.03];
\draw [fill=blue] (7/4,5/2) circle [radius=0.03];
\draw [fill=blue] (5/4,11/4) circle [radius=0.03];
\draw [fill=blue] (3/2,11/4) circle [radius=0.03];

\draw [fill=green] (9/4,9/4) circle [radius=0.03];
\draw [fill=blue] (5/2,9/4) circle [radius=0.03];
\draw [fill=blue] (11/4,9/4) circle [radius=0.03];
\draw [fill=blue] (9/4,5/2) circle [radius=0.03];
\draw [fill=orange] (5/2,5/2) circle [radius=0.03];
\draw [fill=orange] (11/4,5/2) circle [radius=0.03];
\draw [fill=blue] (9/4,11/4) circle [radius=0.03];
\draw [fill=orange] (5/2,11/4) circle [radius=0.03];
\draw [fill=orange] (11/4,11/4) circle [radius=0.03];

\draw [fill=magenta] (13/4,9/4) circle [radius=0.03];
\draw [fill=magenta] (7/2,5/2) circle [radius=0.03];
\draw [fill=orange] (15/4,11/4) circle [radius=0.03];


\draw [fill=blue] (1/2,7/2) circle [radius=0.03];

\draw [fill=orange] (5/4,13/4) circle [radius=0.03];
\draw [fill=blue] (3/2,7/2) circle [radius=0.03];
\draw [fill=orange] (7/4,15/4) circle [radius=0.03];

\draw [fill=magenta] (9/4,13/4) circle [radius=0.03];
\draw [fill=magenta] (5/2,7/2) circle [radius=0.03];
\draw [fill=orange] (11/4,15/4) circle [radius=0.03];

\draw [fill] (7/2,7/2) circle [radius=0.03];

\draw [fill=blue] (5/4,9/4) circle [radius=0.03];
\draw [fill=blue] (9/4,5/4) circle [radius=0.03];
\draw [fill=blue] (11/4,7/4) circle [radius=0.03];
\draw [fill=blue] (7/4,11/4) circle [radius=0.03];

\end{tikzpicture}

\bigskip
c)\quad
\begin{tikzpicture}[scale=1.5]

\draw[help lines,black] (0,0) grid (4,4);

\draw [line width=1.5,->,red] (3/2,1/2) -- (5/2-0.03,1/2);
\draw [line width=1.5,->,red] (1/2,3/2) -- (1/2,5/2-0.03);
\draw [line width=1.5,->,red] (1/2,3/2) -- (3/2-0.03,3/2);
\draw [line width=1.5,->,red] (3/2,1/2) -- (3/2,3/2-0.03);

\draw [line width=1.5,->,green] (1/4,9/4) -- (5/4-0.03,9/4);
\draw [line width=1.5,->,green] (3/4,11/4) -- (3/2-0.03,11/4);
\draw [line width=1.5,->,green] (7/4,5/4) -- (9/4-0.03,5/4);
\draw [line width=1.5,->,green] (7/4,3/2) -- (11/4-0.03,3/2);
\draw [line width=1.5,->,green] (5/4,7/4) -- (5/4,9/4-0.03);
\draw [line width=1.5,->,green] (9/4,1/4) -- (9/4,5/4-0.03);
\draw [line width=1.5,->,green] (11/4,3/4) -- (11/4,3/2-0.03);
\draw [line width=1.5,->,green] (3/2,7/4) -- (3/2,11/4-0.03);
\draw [line width=1.5,->,green] (7/4,7/4) -- (9/4-0.03,7/4);
\draw [line width=1.5,->,green] (7/4,7/4) -- (7/4,9/4-0.03);
\draw [line width=1.5,->,green] (9/4,7/4) -- (9/4,9/4-0.03);
\draw [line width=1.5,->,green] (7/4,9/4) -- (9/4-0.03,9/4);

\draw [line width=1.5,->,blue] (5/4,5/2) -- (9/4-0.03,5/2);
\draw [line width=1.5,->,blue] (11/4,7/4) -- (15/4-0.03,7/4);
\draw [line width=1.5,->,blue] (5/2,5/4) -- (13/4-0.03,5/4);
\draw [line width=1.5,->,blue] (7/4,11/4) -- (9/4-0.03,11/4);
\draw [line width=1.5,->,blue] (5/2,5/4) -- (5/2,9/4-0.03);
\draw [line width=1.5,->,blue] (7/4,11/4) -- (7/4,15/4-0.03);
\draw [line width=1.5,->,blue] (5/4,5/2) -- (5/4,13/4-0.03);
\draw [line width=1.5,->,blue] (11/4,7/4) -- (11/4,9/4-0.03);

\draw [line width=1.5,->,orange] (7/2,3/2) -- (7/2,5/2-0.03);
\draw [line width=1.5,->,orange] (3/2,7/2) -- (5/2-0.03,7/2);
\draw [line width=1.5,->,orange] (5/2,5/2) -- (7/2-0.03,5/2);
\draw [line width=1.5,->,orange] (5/2,5/2) -- (5/2,7/2-0.03);


\draw [fill] (1/2,1/2) circle [radius=0.03];

\draw [fill=red] (5/4,1/4) circle [radius=0.03];
\draw [fill=red] (3/2,1/2) circle [radius=0.03];
\draw [fill=red] (7/4,3/4) circle [radius=0.03];

\draw [fill=green] (9/4,1/4) circle [radius=0.03];
\draw [fill=red] (5/2,1/2) circle [radius=0.03];
\draw [fill=green] (11/4,3/4) circle [radius=0.03];

\draw [fill=blue] (7/2,1/2) circle [radius=0.03];


\draw [fill=red] (1/4,5/4) circle [radius=0.03];
\draw [fill=red] (1/2,3/2) circle [radius=0.03];
\draw [fill=red] (3/4,7/4) circle [radius=0.03];

\draw [fill=green] (5/4,5/4) circle [radius=0.03];
\draw [fill=green] (7/4,5/4) circle [radius=0.03];
\draw [fill=red] (3/2,3/2) circle [radius=0.03];
\draw [fill=green] (7/4,3/2) circle [radius=0.03];
\draw [fill=green] (5/4,7/4) circle [radius=0.03];
\draw [fill=green] (3/2,7/4) circle [radius=0.03];
\draw [fill=green] (7/4,7/4) circle [radius=0.03];
\draw [fill=blue] (11/4,7/4) circle [radius=0.03];

\draw [fill=blue] (5/2,5/4) circle [radius=0.03];
\draw [fill=green] (11/4,3/2) circle [radius=0.03];
\draw [fill=green] (9/4,7/4) circle [radius=0.03];

\draw [fill=blue] (13/4,5/4) circle [radius=0.03];
\draw [fill=orange] (7/2,3/2) circle [radius=0.03];
\draw [fill=blue] (15/4,7/4) circle [radius=0.03];


\draw [fill=green] (1/4,9/4) circle [radius=0.03];
\draw [fill=red] (1/2,5/2) circle [radius=0.03];
\draw [fill=green] (3/4,11/4) circle [radius=0.03];

\draw [fill=green] (7/4,9/4) circle [radius=0.03];
\draw [fill=blue] (5/4,5/2) circle [radius=0.03];
\draw [fill=green] (3/2,11/4) circle [radius=0.03];
\draw [fill=blue] (7/4,11/4) circle [radius=0.03];

\draw [fill=green] (9/4,9/4) circle [radius=0.03];
\draw [fill=blue] (5/2,9/4) circle [radius=0.03];
\draw [fill=blue] (11/4,9/4) circle [radius=0.03];
\draw [fill=blue] (9/4,5/2) circle [radius=0.03];
\draw [fill=orange] (5/2,5/2) circle [radius=0.03];
\draw [fill=blue] (9/4,11/4) circle [radius=0.03];
\draw [fill=orange] (11/4,11/4) circle [radius=0.03];

\draw [fill=magenta] (13/4,9/4) circle [radius=0.03];
\draw [fill=orange] (7/2,5/2) circle [radius=0.03];
\draw [fill=magenta] (15/4,11/4) circle [radius=0.03];


\draw [fill=blue] (1/2,7/2) circle [radius=0.03];

\draw [fill=blue] (5/4,13/4) circle [radius=0.03];
\draw [fill=orange] (3/2,7/2) circle [radius=0.03];
\draw [fill=blue] (7/4,15/4) circle [radius=0.03];

\draw [fill=magenta] (9/4,13/4) circle [radius=0.03];
\draw [fill=orange] (5/2,7/2) circle [radius=0.03];
\draw [fill=magenta] (11/4,15/4) circle [radius=0.03];

\draw [fill] (7/2,7/2) circle [radius=0.03];


\draw [fill=green] (3/2,5/4) circle [radius=0.03];
\draw [fill=green] (5/4,3/2) circle [radius=0.03];
\draw [fill=green] (5/4,9/4) circle [radius=0.03];
\draw [fill=green] (9/4,5/4) circle [radius=0.03];
\draw [fill=blue] (11/4,5/4) circle [radius=0.03];
\draw [fill=blue] (9/4,3/2) circle [radius=0.03];
\draw [fill=blue] (5/2,3/2) circle [radius=0.03];
\draw [fill=blue] (5/2,7/4) circle [radius=0.03];
\draw [fill=blue] (3/2,9/4) circle [radius=0.03];
\draw [fill=blue] (3/2,5/2) circle [radius=0.03];
\draw [fill=blue] (7/4,5/2) circle [radius=0.03];
\draw [fill=blue] (5/4,11/4) circle [radius=0.03];
\draw [fill=orange] (11/4,5/2) circle [radius=0.03];
\draw [fill=orange] (5/2,11/4) circle [radius=0.03];

\end{tikzpicture}\qquad d)\quad
\begin{tikzpicture}[scale=1.5]

\draw[help lines,black] (0,0) grid (4,4);

\draw [line width=1.5,->,red] (3/2,1/2) -- (5/2-0.03,1/2);
\draw [line width=1.5,->,red] (1/2,3/2) -- (1/2,5/2-0.03);
\draw [line width=1.5,->,red] (1/2,3/2) -- (5/4-0.03,3/2);
\draw [line width=1.5,->,red] (3/2,1/2) -- (3/2,5/4-0.03);

\draw [line width=1.5,->,green] (1/4,9/4) -- (5/4-0.03,9/4);
\draw [line width=1.5,->,green] (3/4,11/4) -- (3/2-0.03,11/4);
\draw [line width=1.5,->,green] (7/4,5/4) -- (9/4-0.03,5/4);
\draw [line width=1.5,->,green] (7/4,3/2) -- (11/4-0.03,3/2);
\draw [line width=1.5,->,green] (5/4,7/4) -- (5/4,9/4-0.03);
\draw [line width=1.5,->,green] (9/4,1/4) -- (9/4,5/4-0.03);
\draw [line width=1.5,->,green] (11/4,3/4) -- (11/4,3/2-0.03);
\draw [line width=1.5,->,green] (3/2,7/4) -- (3/2,11/4-0.03);
\draw [line width=1.5,->,green] (7/4,7/4) -- (9/4-0.03,7/4);
\draw [line width=1.5,->,green] (7/4,7/4) -- (7/4,9/4-0.03);
\draw [line width=1.5,->,green] (9/4,7/4) -- (9/4,9/4-0.03);
\draw [line width=1.5,->,green] (7/4,9/4) -- (9/4-0.03,9/4);

\draw [line width=1.5,->,blue] (5/4,5/2) -- (9/4-0.03,5/2);
\draw [line width=1.5,->,blue] (11/4,7/4) -- (15/4-0.03,7/4);
\draw [line width=1.5,->,blue] (5/2,5/4) -- (13/4-0.03,5/4);
\draw [line width=1.5,->,blue] (7/4,11/4) -- (9/4-0.03,11/4);
\draw [line width=1.5,->,blue] (5/2,5/4) -- (5/2,9/4-0.03);
\draw [line width=1.5,->,blue] (7/4,11/4) -- (7/4,15/4-0.03);
\draw [line width=1.5,->,blue] (5/4,5/2) -- (5/4,13/4-0.03);
\draw [line width=1.5,->,blue] (11/4,7/4) -- (11/4,9/4-0.03);

\draw [line width=1.5,->,orange] (7/2,3/2) -- (7/2,5/2-0.03);
\draw [line width=1.5,->,orange] (3/2,7/2) -- (5/2-0.03,7/2);
\draw [line width=1.5,->,orange] (11/4,5/2) -- (7/2-0.03,5/2);
\draw [line width=1.5,->,orange] (5/2,11/4) -- (5/2,7/2-0.03);


\draw [fill] (1/2,1/2) circle [radius=0.03];

\draw [fill=red] (5/4,1/4) circle [radius=0.03];
\draw [fill=red] (3/2,1/2) circle [radius=0.03];
\draw [fill=red] (7/4,3/4) circle [radius=0.03];

\draw [fill=green] (9/4,1/4) circle [radius=0.03];
\draw [fill=red] (5/2,1/2) circle [radius=0.03];
\draw [fill=green] (11/4,3/4) circle [radius=0.03];

\draw [fill=blue] (7/2,1/2) circle [radius=0.03];


\draw [fill=red] (1/4,5/4) circle [radius=0.03];
\draw [fill=red] (1/2,3/2) circle [radius=0.03];
\draw [fill=red] (3/4,7/4) circle [radius=0.03];

\draw [fill=green] (5/4,5/4) circle [radius=0.03];
\draw [fill=red] (3/2,5/4) circle [radius=0.03];
\draw [fill=green] (7/4,5/4) circle [radius=0.03];
\draw [fill=red] (5/4,3/2) circle [radius=0.03];
\draw [fill=green] (3/2,3/2) circle [radius=0.03];
\draw [fill=green] (7/4,3/2) circle [radius=0.03];
\draw [fill=green] (5/4,7/4) circle [radius=0.03];
\draw [fill=green] (3/2,7/4) circle [radius=0.03];
\draw [fill=green] (7/4,7/4) circle [radius=0.03];
\draw [fill=blue] (11/4,7/4) circle [radius=0.03];

\draw [fill=blue] (5/2,5/4) circle [radius=0.03];
\draw [fill=green] (11/4,3/2) circle [radius=0.03];
\draw [fill=green] (9/4,7/4) circle [radius=0.03];

\draw [fill=blue] (13/4,5/4) circle [radius=0.03];
\draw [fill=orange] (7/2,3/2) circle [radius=0.03];
\draw [fill=blue] (15/4,7/4) circle [radius=0.03];


\draw [fill=green] (1/4,9/4) circle [radius=0.03];
\draw [fill=red] (1/2,5/2) circle [radius=0.03];
\draw [fill=green] (3/4,11/4) circle [radius=0.03];

\draw [fill=green] (7/4,9/4) circle [radius=0.03];
\draw [fill=blue] (5/4,5/2) circle [radius=0.03];
\draw [fill=green] (3/2,11/4) circle [radius=0.03];
\draw [fill=blue] (7/4,11/4) circle [radius=0.03];

\draw [fill=green] (9/4,9/4) circle [radius=0.03];
\draw [fill=blue] (5/2,9/4) circle [radius=0.03];
\draw [fill=blue] (11/4,9/4) circle [radius=0.03];
\draw [fill=blue] (9/4,5/2) circle [radius=0.03];
\draw [fill=orange] (5/2,5/2) circle [radius=0.03];
\draw [fill=orange] (11/4,5/2) circle [radius=0.03];
\draw [fill=blue] (9/4,11/4) circle [radius=0.03];
\draw [fill=orange] (5/2,11/4) circle [radius=0.03];
\draw [fill=orange] (11/4,11/4) circle [radius=0.03];

\draw [fill=magenta] (13/4,9/4) circle [radius=0.03];
\draw [fill=orange] (7/2,5/2) circle [radius=0.03];
\draw [fill=magenta] (15/4,11/4) circle [radius=0.03];


\draw [fill=blue] (1/2,7/2) circle [radius=0.03];

\draw [fill=blue] (5/4,13/4) circle [radius=0.03];
\draw [fill=orange] (3/2,7/2) circle [radius=0.03];
\draw [fill=blue] (7/4,15/4) circle [radius=0.03];

\draw [fill=magenta] (9/4,13/4) circle [radius=0.03];
\draw [fill=orange] (5/2,7/2) circle [radius=0.03];
\draw [fill=magenta] (11/4,15/4) circle [radius=0.03];

\draw [fill] (7/2,7/2) circle [radius=0.03];

\draw [fill=green] (5/4,9/4) circle [radius=0.03];
\draw [fill=green] (9/4,5/4) circle [radius=0.03];
\draw [fill=blue] (11/4,5/4) circle [radius=0.03];
\draw [fill=blue] (9/4,3/2) circle [radius=0.03];
\draw [fill=blue] (5/2,3/2) circle [radius=0.03];
\draw [fill=blue] (5/2,7/4) circle [radius=0.03];
\draw [fill=blue] (3/2,9/4) circle [radius=0.03];
\draw [fill=blue] (3/2,5/2) circle [radius=0.03];
\draw [fill=blue] (7/4,5/2) circle [radius=0.03];
\draw [fill=blue] (5/4,11/4) circle [radius=0.03];

\end{tikzpicture}

\end{center}
\caption{\label{fig1} Representation of the structure of the complexes of complex-valued differential forms related to the Iwasawa manifold and some of its small deformations. Horizontal and vertical arrows represent respectively the maps $\p$ and $\pb$. The different colours correspond to the various steps of the Stelzig process.}
\end{figure}
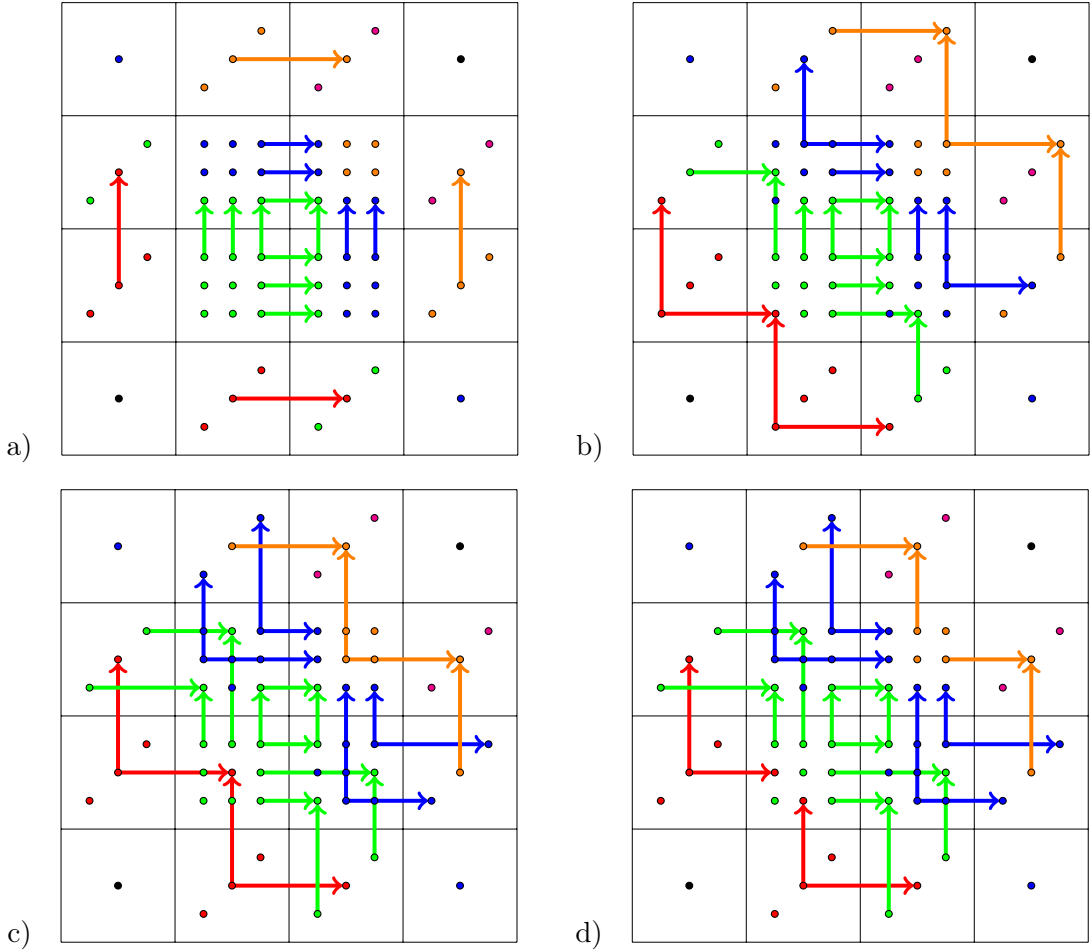

In Figure \ref{fig1}, we see the graphs corresponding to some deformations of the Iwasawa manifold obtained in this way. The graph a) corresponds to the complex of differential forms of the Iwasawa manifold itself. Its structure, already determined in \cite[Section 4]{Ang13} by D. Angella is quite simple and it can be represented as direct sum of thirtysix zigzags of lenght $1$, twelve zigzags of lenght $2$ and a square.
The graph b) instead represents the complex related to a deformation of the Iwasawa manifold in the case where 
\[
t_{21}\neq 0,\quad t_{11}=t_{12}=t_{22}=t_{31}=t_{32}=0,
\]  
that is of class (ii.a). It has a bit more complicated structure, with twentyeight zigzags of lenght $1$, four zigzags of lenght $2$,  four zigzags of lenght $3$, two zigzags of lenght $5$ and a square.

Finally the graphs c) and d) concern a deformation obtained by parameters
\[
t_{11}\neq 0,\quad t_{22}\neq0,\quad t_{12}=t_{21}=t_{31}=t_{32}=0
\]
respectively when $|t_{11}|=|t_{22}|$ and $|t_{11}|\neq|t_{22}|$. In the first case the deformation represents an example of class (iii.a), while in the second case it belongs to the class (iii.b).
The only difference beetween the two graphs is given by two zigzags of lenght $5$ in the first, which split into four zigzags of lenght $3$ if the two norms are different. In addition to these, we can find also eight more zigzags of lenght $3$ and a square.

Because of the results mentioned before, these representations give another way to compute the dimensions of de Rham, Dolbeault, Bott-Chern and Aeppli cohomologies of the manifolds considered, already determined in \cite[p. 96]{Nak75} and \cite[Section 5]{Ang13}. For example, if we want to check the dimensions of Dolbeault cohomology, we have to remove vertical arrows in the graphs and count the remaining points (see the table in \cite[p. 96]{Nak75} to check the corresponding values).

\section{Properties of the Fr\"olicher spectral sequence under
 small deformations}
Thanks to the explicit description given in \cite{CFUG97}, we can easily use the graphs in Figure \ref{fig1} to determine the dimensions of the successive pages of the Fr\"olicher spectral sequence of the related manifolds, that we collect in the table of Figure \ref{fig2}.
In particular, we observe that for the Iwasawa Manifold and the deformation of class (ii) analysed the Fr\"olicher spectral sequence degenerates at the second page, while in the case of the deformations
 of class (iii), in both the examples considered, it degenerates at the first page, according to the Hodge numbers of the Iwasawa manifold and its small deformations already determined by I. Nakamura in \cite[p.96]{Nak75}.
These facts permit us to make some considerations. 

\begin{re}{1}
M. Maschio showed in \cite[Theorem 1]{Mas18} that, given any compact complex manifold $X$, if the dimensions of Dolbeault cohomology are constant under deformations and its Fr\"olicher spectral sequence degenerates at the second page, then for each small enough deformation of $X$ its Fr\"olicher spectral sequence degenerates at the second page too. In \cite[Section 5]{Mas18} he also provided that the hypothesis of constantness is strictly necessary to guarantee  the stability of degeneration at the second page under small deformations: indeed, using the Nakamura manifold, that represents an example of holomorphic parallelizable complex three-dimensional solvmanifold and whose Fr\"olicher spectral sequence degenerates at the second page, he considered a family of small deformations and proved  (\cite[Theorem 9]{Mas18}) that their Fr\"olicher spectral sequence degenerates at higher steps.

Here we provided, by the Iwasawa manifold and its small deformations, another example to show 
the necessity of this condition,
using in this case a nilmanifold with a quite simple structure.
\end{re}

\begin{figure}
{\small
\[
\begin{array}{|c|cc|ccc|cccc|ccc|cc|}
\toprule
\text{a)} & {1,0} & {0,1} &{2,0} & {1,1} &{0,2} & {3,0} &{2,1} & {1,2} & {0,3} & {3,1} &{2,2} & {1,3} &{3,2} & {2,3} \\
\midrule
E_1^{\bullet,\bullet} & 3 & 2 & 3 & 6 & 2 & 1 & 6 & 6 & 1 & 2 & 6 & 3 & 2 & 3\\
E_2^{\bullet,\bullet} & 2 & 2 & 2 & 4 & 2 & 1 & 4 & 4 & 1 & 2 & 4 & 2 & 2 & 2\\
\bottomrule
\end{array}
\]
\[
\begin{array}{|c|cc|ccc|cccc|ccc|cc|}
\toprule
\text{b)} & {1,0} & {0,1} &{2,0} & {1,1} &{0,2} & {3,0} &{2,1} & {1,2} & {0,3} & {3,1} &{2,2} & {1,3} &{3,2} & {2,3}\\
\midrule
E_1^{\bullet,\bullet} & 2 & 2 & 2 & 5 & 2 & 1 & 5 & 5 & 1 & 2 & 5 & 2 & 2 & 2\\
E_2^{\bullet,\bullet} & 2 & 2 & 2 & 4 & 2 & 1 & 4 & 4 & 1 & 2 & 4 & 2 & 2 & 2
\\
\bottomrule
\end{array}
\]
\[
\begin{array}{|c|cc|ccc|cccc|ccc|cc|}
\toprule
\text{c), d)} & {1,0} & {0,1} &{2,0} & {1,1} &{0,2} & {3,0} &{2,1} & {1,2} & {0,3} & {3,1} &{2,2} & {1,3} &{3,2} & {2,3} \\
\midrule
E_1^{\bullet,\bullet} & 2 & 2 & 1 & 5 & 2 & 1 & 4 & 4 & 1 & 2 & 5 & 1 & 2 & 2\\
\bottomrule
\end{array}
\]
}
\caption{\label{fig2} Successive pages of the Fr\"olicher spectral sequence of the Iwasawa manifold (table a)) and some of its small deformations (tables b) and c)).}
\end{figure}

\begin{re}{2}
The example of the Nakamura manifold used by M. Maschio permits also to observe (\cite[Remark 2]{Mas18}) that, in general, the dimensions of the components $E_r^{p,q}(X_\tg)$ of the Fr\"olicher spectral sequence do not vary either upper
semi-continuously or lower semi-continuously under small deformations, i.e. if $\{X_{\mathbf{t}}=(X,J_{\mathbf{t}})\}_{\mathbf{t}\in\Delta(0,\varepsilon)}$ is a family of small deformations of a compact complex manifold $X$, then for every $p,q\in\Z$ and $r\in\Z\setminus\{0,1\}$ the map 
\[
\tg\mapsto \dim\big(E_r^{p,q}(X_\tg)\big)
\]
is, in general, nor upper neither lower semi-continuos.

We can deduce it also by analysing values exposed in Figure \ref{fig2}, in particular comparing the numbers corresponding to the components of type $(2,0)$ and $(1,1)$ of tables a) and c). 
\end{re}

\begin{re}{3}
It is evident that the Fr\"olicher spectral sequences of the considered manifolds present also a symmetry property: for each $p,q\in \N$ we have 
\[
\dim\big(E_2^{p,q}(X_\tg)\big)=\dim\big(E_2^{3-p,3-q}(X_\tg)\big).
\]
This equality represents a generalization of the Serre duality for Dolbeault cohomology, correspponding for each $p,q\in\Z$ to the isomorphism
\[
H_{\pb}^{p,q}(X_\tg)\simeq H_{\pb}^{n-p,n-q}(X_\tg)^*.
\]
It was already noticed by L. Ugarte in \cite{Uga00} for a hypothetical complex structure on $S^6$ and it has been recently proved by A. Milojevic in \cite[Corollary 7]{Mil19}, where he shows that Serre symmetry is valid in general for every page of the Fr\"olicher spectral sequence.
\end{re}

\end{document}